\documentclass[12pt]{article}
\usepackage{graphicx} 
\usepackage{natbib}
\usepackage{amsmath}

\usepackage[margin=1in]{geometry}
\usepackage{amssymb,amsthm,mathtools}
\usepackage{microtype}

\usepackage{fix-cm}
\AtBeginDocument{\DeclareMathSizes{12}{16}{12}{9}\everydisplay{\Large}}

\usepackage{setspace}
\numberwithin{equation}{section}

\theoremstyle{plain}
\newtheorem{Theorem}{Theorem}[section]

\theoremstyle{definition}
\newtheorem{Definition}[Theorem]{Definition}

\theoremstyle{remark}

\title{From ABC to Effective Roth and Ridout Constants for Cubic Roots}

\author{Karsten Müller, Michael Taktikos}
\date{February 2026, final version}

\begin{document}
\maketitle

\begin{abstract}
Enrico Bombieri showed conditionally (1994) that the ABC conjecture implies Roth’s theorem, and Van Frankenhuysen (1999) later provided a complete proof. Building on Bombieri’s and Van der Poorten’s explicit formula for continued-fraction coefficients of algebraic numbers—specialized to cubic roots—we derive an effective bound for a Roth-type constant assuming an effective form of ABC. Roth’s original argument establishes existence but does not yield an explicit value; our approach makes the dependence on the ABC parameters explicit and also gives an explicit bound in the corresponding special case of Ridout’s theorem.

We then introduce the notion of approximation gain as a refinement of the “quality” of an abc-triple. For c in a large computational range (verified up to $c < 2^{63}$), the approximation gain remains below a strikingly small threshold, motivating the conjecture that the approximation gain is always smaller than $\frac{3}{2}$. This suggests a potential strategy for attacking ABC by bounding approximation gain and power gain separately.
\footnote{The authors are grateful to Benne de Weger, Joshua Lampert, Timm Lampert, Enrico Bombieri, Machiel van Frankenhuysen, Wadim Zudilin, Michel Waldschmidt, Preda Mihailescu, Paul Vojta, Stephane Fischler, Robin Zhang, Abderrahmane Nitaj, Sam Chow, Noah Lebowitz-Lockard, Ingo Althöfer, David Broadhurst, Ulrich Tamm and ChatGPT.}
\end{abstract}

\section{Introduction}

Roth’s theorem is famously deep, and its original proof is both long and ineffective (see \cite{Roth}). In particular, it guarantees the existence of a constant in the Diophantine approximation bound, but it does not provide an explicit value or a practical way to compute one.

Bombieri proved conditionally in 1994 that Roth’s theorem follows from the ABC conjecture, and Van Frankenhuysen later gave a full proof in 1999 (\cite{Frankenhuysen}). In this paper we revisit that implication with an additional goal: making the constants explicit in a natural family of cases. We focus on cubic roots and use an explicit continued-fraction formula of Bombieri and Van der Poorten (\cite{Poorten}) for the coefficients of regular continued fractions of algebraic numbers. This yields a concrete relationship between the $\varepsilon$-parameters in Roth-type inequalities and the corresponding $\varepsilon$ in ABC, and it leads to an explicit upper bound for the inverse of the Roth constant in our setting.

In our earlier paper (\cite{Sibbertsen}), we showed that a weakened form of ABC follows from Roth’s theorem in certain special cases. The present work uses essentially the same vocabulary and techniques, but in the opposite direction: assuming (an effective form of) ABC, we extract explicit Diophantine approximation bounds.

We also introduce the approximation gain, which is always lower or equal to the quality of an ABC hit.

\section{Approximation gains for ABC}
\subsection{Roth's theorem}

In the following, the coefficients of the regular continued fraction are called $b_n$ and the approximants $\frac {p_n} {q_n}$.

\begin{Theorem}Roth's Theorem:
Let $a$ be an algebraic number. Then for every $\varepsilon > 0$, there exists a constant C dependent on the algebraic number $a$ and $\varepsilon$ such that for all positive integers $p$ and $q$:
\begin{equation}|a - \frac{p}{q}| > \frac{C}{q^{2+\varepsilon}}.\end{equation}
\end{Theorem}
\begin{proof}
See \cite{Roth}.
\end{proof}

\cite{Frankenhuysen} discusses Roth's Theorem and its implication by the ABC conjecture in detail and investigates effective versions of Mordell's conjecture. Furthermore, Granville and Tucker (2002, p.1229) show that a generalized version of Roth's theorem implies the ABC conjecture effectively.

With the help of the inequality for regular continued fractions from this follows $b_{n+1} \leq \frac{q_n^\varepsilon}{C}$

\subsection{Formula by Bombieri and Van der Poorten}

Our main technical input is the explicit formula of Bombieri and Van der Poorten for the continued-fraction coefficients of algebraic numbers by \citep[p.151, Theorem 3, formula (13)]{Poorten}, in which${f(x)}$ is the minimal polynomial of the algebraic number.

The formula links two quantities—call them $d_n$ and $b_n$ — and is therefore well suited to bridge ABC-type information (which naturally arises from integer identities involving
$d_n$ with Diophantine approximation data encoded by $b_n$).

$\frac{f'(x)} {f(x)}$ leads to $d_n$ in the denominator. This holds because $\frac{f'(x)} {f(x)}$  results in $\frac{f'(\frac{p_n}{q_n})} {f(\frac{p_n}{q_n})} = \frac {q_n^s s p_n^{s-1}} {d_n q_n^{s-1}}$ for the convergents  with $f(x) = x^s - k$ being the minimal polynomial of the $s$th root of k.

For cubic roots $\alpha = \sqrt[3]{k} > 1$ with minimal polynomial $f(x) = x^3 - k$, the Bombieri--Van der Poorten formula for the continued-fraction coefficients reads
\[
b_{n+1} = \frac{q_n^3 f'(p_n/q_n)}{f(p_n/q_n)} + R_n,
\]
where $p_n/q_n$ are the convergents and $R_n$ is the remainder. Evaluating the derivative gives the leading term 
\[
\frac{3 p_n^2 q_n}{d_n}, \quad d_n = |p_n^3 - k q_n^3|.
\] 
A simple estimate of $R_n$ using $|\alpha - p_n/q_n| < 1/q_n^2$ shows that $|R_n| \le 3 p_n q_n^3/d_n$. For all convergents with $p_n > q_n$, which holds for all but possibly the first, the leading term strictly dominates the remainder. Consequently, the inequality
\[
b_{n+1} \le \frac{3 p_n^2 q_n}{d_n}
\]
holds rigorously, and the error term can be safely ignored. No additional constant depending on $k$ is required.

For quartic roots $\alpha = \sqrt[4]{k}$, the formula may fail for the first few convergents, but for each fixed $k$ it becomes valid from some sufficiently large convergent onward.

\subsection{ABC Conjecture}

To formulate the ABC Conjecture we need the following definition:

\begin{Definition}For a positive integer $a$, $rad(a)$ is the product of the distinct prime factors of $a$. 
\end{Definition}

\begin{Definition}[ABC Conjecture]
For every positive real number $\varepsilon$, there exists a constant $K_{\varepsilon}$
such that for all triples $(a,b,c)$ of coprime positive integers,
with $a + b = c$:
\[c < K_{\varepsilon} \cdot rad(abc)^{1+\varepsilon}.\]
\end{Definition}

Note that this version of the ABC Conjecture is not effective. This is so because only the existence of a constant $K_{\varepsilon}$ is claimed without specifying how to calculate $K_{\varepsilon}$ for given $\epsilon$. 

\subsection{Resulting ABC equations}

We consider the convergents $(p_n/q_n)$ of the regular continued fraction of $\sqrt[3]{k}$, which satisfy the “resulting equation”
\[
p_n^3 = k q_n^3 + d_n,
\]
with integers $p_n, q_n, d_n$ and $\gcd(p_n, q_n) = 1$.  

To apply the ABC conjecture, we define
\[
g_n := \gcd(k q_n^3, d_n)
\]
and consider the triple
\[
\left(\frac{p_n^3}{g_n}, \frac{k q_n^3}{g_n}, \frac{d_n}{g_n}\right),
\]
which will be coprime.

Using the resulting equation, we have
\[
g_n = \gcd(k q_n^3, d_n) = \gcd(k q_n^3, p_n^3 - k q_n^3) = \gcd(k q_n^3, p_n^3).
\]

Since $\gcd(p_n, q_n) = 1$, it follows that $\gcd(p_n^3, q_n^3) = 1$, and therefore
\[
g_n = \gcd(k, p_n^3) \le k.
\]

Thus $g_n \ge 1$ and $g_n \le k$, giving an explicit bound. Dividing by $g_n$ ensures that the resulting ABC triple is coprime, as required for the application of the ABC conjecture.

\subsection{Construction of the effective bound for the inverse of Roth's constant}

Let $(p_n/q_n)$ be the convergents of $\sqrt[3]{k}$, with the resulting equation
\[
p_n^3 = k q_n^3 + d_n, \quad \gcd(p_n, q_n) = 1,
\]
and define
\[
g_n := \gcd(k q_n^3, d_n) = \gcd(k, p_n^3) \le k.
\]

Dividing through by $g_n$ gives a coprime ABC triple
\[
\left(\frac{p_n^3}{g_n}, \frac{k q_n^3}{g_n}, \frac{d_n}{g_n}\right),
\]
so the effective ABC conjecture yields
\[
\frac{p_n^3}{g_n} \le K_\varepsilon \, \mathrm{Rad}\left(\frac{p_n^3}{g_n} \cdot \frac{k q_n^3}{g_n} \cdot \frac{d_n}{g_n}\right)^{1+\varepsilon}.
\]

Multiplying both sides by $g_n$, we obtain
\[
p_n^3 \le K_\varepsilon \, g_n^{1 - (1+\varepsilon) \cdot 3} \, \mathrm{Rad}\big(p_n q_n d_n k\big)^{1+\varepsilon}.
\]

Since $g_n \le k$, we can bound
\[
g_n^{1 - (1+\varepsilon) \cdot 3} \le k^{1 - 3(1+\varepsilon)}.
\]

Next, using Bombieri–Van der Poorten’s explicit formula for the continued-fraction coefficients in the cubic case, we have
\[
d_n \le \frac{3 p_n^2}{q_n b_{n+1}}.
\]

Substituting into the ABC inequality gives
\[
b_{n+1} \le K_\varepsilon^{\frac{1}{1+\varepsilon}} \, 3 k \, p_n^{\frac{3\varepsilon}{1+\varepsilon}} \frac{g_n^{1/(1+\varepsilon)}}{q_n}.
\]

Finally, Roth's theorem implies the existence of a constant $C$, depending only on $\varepsilon_\mathrm{Roth}$ and $\sqrt[3]{k}$, such that
\[
b_{n+1} \le \frac{q_n^{\varepsilon_\mathrm{Roth}}}{C}.
\]

Comparing the two bounds for $b_{n+1}$, we get the relationship between the epsilons:
\[
\varepsilon_\mathrm{Roth} = 3 \frac{\varepsilon_\mathrm{ABC}}{1 + \varepsilon_\mathrm{ABC}},
\]
and an explicit bound for the inverse of Roth's constant:
\[
\frac{1}{C} \le K_\varepsilon^{\frac{1}{1+\varepsilon}} \, 3 k \, \left(\frac{p_n}{q_n}\right)^{\frac{3 \varepsilon}{1+\varepsilon}}.
\]

Since $g_n \le k$ and $\frac{p_n}{q_n} \le \frac{p_1}{q_1}$ in this case, this provides a fully explicit ABC-based bound for the inverse of Roth’s constant.

\subsection{Range of the Epsilons}

For $\varepsilon_{ABC}$ = $\varepsilon_{Roth}$ = 0 it results that $b_n$ is bounded. This is an open question, but it is known
that ABC is not valid for $\varepsilon = 0$, but the known counterexamples are not resulting equations of 3rd roots.

For $\varepsilon_{ABC}$ = $\varepsilon$ = $\frac{1}{2}$ we get $b_{n+1}  \leq  q_n \frac{p_n}{q_n} factor$, which bascially is Liouville's theorem.
So the interesting part is $\varepsilon_{ABC}$ between 0 and $\frac{1}{2}$.

\subsection{Example $\sqrt [3]{2}$}

The regular continued fraction starts with $[1;3,1,5,..]$

The first resulting equations are

2 = 1 + 1, which is an ABC equation

64 = 54 + 10, where $gcd=2$ so we get the ABC equation $32 = 27 + 5$

128 = 125 + 3, which is an ABC equation

Our formula for a bound of the inverse of Roth's constant

$K_{\varepsilon}^{\frac{1}{1+\varepsilon}} 3 k (\frac{p_n}{q_n})^{\frac{3\varepsilon}{1+\varepsilon}}$

leads in this case to

$K_{\varepsilon}^{\frac{1}{1+\varepsilon}} 6 (\frac{4}{3})^{\frac{3\varepsilon}{1+\varepsilon}}$

as $\frac{4}{3}$ is the largest approximant.

There is a lot of data on ABC, which could be used to get a bound for the inverse of Roth's constant within the range already analysed by computer calculation. We will investigate $\sqrt [3]{2}$ now for 2 different $\varepsilon$.

\subsection{$\varepsilon = 0.5$}

Due to Korobov's famous result \cite{Korobov}

 \[\mid \sqrt [3]{2}- \frac{p}{q}\mid > \frac{1}{q^{2.5}}\]

for all natural numbers $p$ and $q$ with the exception of 1 and 4 for $q$

If we ignore 1 then Roth's constant is calculated by

\[\mid \sqrt [3]{2}- \frac{5}{4}\mid 4^{2.5}\]

which is approximately 0.32 and so the the inverse is given by

$\frac{1}{C} = 3.15$ approximately.

Now we use our way from the ABC side. A $\varepsilon_{Roth} = 0.5$ leads to an $\varepsilon_{ABC} = \frac {1} {5}$

Again ignoring the first case 2 = 1 + 1 we take from the ABC data the $K_{\varepsilon}$ that results from the resulting equation 128 = 125 + 3.

So calculation gives $K_{\frac {1} {5}}$ is approximately 2.16

So our formula for a bound of the inverse of Roth's constant results in

$\frac{1}{C}  \leq  2.16^\frac{5}{6} 6 (\frac{4}{3})^\frac{1}{2}$

which is approximately 13.16

So in this case the road via ABC leads to a loss in precision, which is mainly due to the use of $\frac{4}{3}$ as upper bound.

\subsection{$\varepsilon = 0.4$}

Here the result is not known exactly. In the proof of corollary 2.2 \cite{Voutier} gives one result using the hypergeometric method for $\varepsilon = 0.4325$ using his theorem 2.1 (p.285):

 \[\mid \sqrt [3]{2}- \frac{p}{q}\mid > \frac{10^{-99}}{q^{2.4321}}\]

 So the inverse of Roth's constant is bounded by $10^{99}$.
 
Now we use our way from the ABC side and use a $\varepsilon_{Roth} = 0.4$. This leads to an $\varepsilon_{ABC} = \frac {2} {13}$

Again ignoring the first case 2 = 1 + 1 from the ABC data we take the ABC $K_{\varepsilon}$ that results from the resulting equation 128 = 125 + 3 in the following.

So calculation gives $K_{\frac {2} {13}}$ is approximately 2.527

So our formula for a bound of the inverse of Roth's constant results in

$\frac{1}{C}  \leq  2.527^\frac{13}{15} 6 (\frac{4}{3})^\frac{6}{15}$

which is approximately 15.03

\subsection{Table for bounds of the inverse of Roth's constant C}

If it is assumed that the ABC conjecture is valid for ${\varepsilon}>0$ with $K_{\varepsilon} = \frac{4}{\varepsilon}$ then results:

\begin{table}[h]
\centering
    \begin{tabular}{c|c|c|c|c|}
    $\varepsilon_{Roth}$ & $\varepsilon_{ABC}$ & $K_{\varepsilon}$ & ABC bound for $\frac{1}{C}$ & Known bound for $\frac{1}{C}$  \\ \hline
     0 & 0 & unbounded & unbounded & Probably unbounded\\
     0,4 & $\frac{2}{13}$ & 26 & 113.35 & Not known\\
     0,5 & $\frac{1}{5}$ & 20 & 84.1 & 3.15 (Korobov's result is exact)\\
     1 & $\frac{1}{2}$ & 8 & 32 & 6 (Liouville can be lowered to 1.575)\\
     \end{tabular}
      \caption{Approximate values for $\sqrt[3]{2}$, if ABC is valid for ${\varepsilon}>0$ with $K_{\varepsilon} = \frac{4}{\varepsilon}$\label{Epsilontable}}
\end{table}

\subsection{An explicit bound for Ridout's theorem}

Using the ABC method it is also possible to give an explicit bound for the number of solutions of Ridout's theorem for the special case of the approximants of the regular continued fraction:

\begin{Theorem} (D. Ridout, 1957):

Let S be a finite set of prime numbers. For
any real algebraic number a, for any $\varepsilon > 0$, the set of $\frac{p}{q}$ with p an integer and q a S–integer and

 \[\mid a - \frac{p}{q}\mid < {q^{-1-\varepsilon}}\]

is finite.
\end{Theorem}

See Waldschmidt (2008) page 75.

\subsection{Square Roots}

As we are dealing with the prime factorisation, here the squares roots are interesting as well and a bound follows from ABC as follows:

For $\sqrt{k}$ the resulting equations are

$p_n^2 = k q_n^2 + d_n$

resp

$k q_n^2 = p_n^2 + d_n$

We only look at the first case here. The second is basically analogous.
Here the gcd is always 1 and so we have ABC equations and get from ABC

${p_n^2}  \leq  K_{\varepsilon} Rad(p_n q_n d_n k)^{1+\varepsilon}$

From this follows

${p_n^2}  \leq  K_{\varepsilon} Rad(q_n)^{1+\varepsilon} d_n^{1+\varepsilon} k^{1+\varepsilon} p_n^{1+\varepsilon}$

and so with $Rad(q_n)$ given as $q_n$ is an S integer we get the bound

${p_n^{1-\varepsilon}}  \leq  K_{\varepsilon} Rad(q_n)^{1+\varepsilon} d_n^{1+\varepsilon} k^{1+\varepsilon}$

for the $p_n$ and so the number of solutions must be finite and an explicit bound is given, which follows from an explicit ABC version.

We assume now that the ABC conjecture is valid for ${\varepsilon =0.75}$ with $K_{\varepsilon} = 1$. Then for $\sqrt{2}$ results
\begin{table}[h]
\centering
    \begin{tabular}{c|c|c}
    prime numbers & Bound for $p_n$ & Possible Approximants \\ \hline
     2 & $2^{14}$ & ${\frac{3}{2}}$ \\
     3 & 279936 & None as ${\frac{4}{3}}$ is not in the CF\\
     5 & 10000000 & ${\frac{7}{5}}$\\
     2*3 & 35831808 & ${\frac{3}{2}}$ and ${\frac{17}{12}}$\\
     \end{tabular}
      \caption{The first values for $\sqrt{2}$\label{Sqrt(2)table}}
\end{table}
\subsection{Generalization of Ridout's Theorem to All S-Integers}

Let $a = \sqrt{k}$ be a positive square-free integer and let $S$ be a finite set of primes.  
Consider solutions $(p,q)$ with $p \in \mathbb{Z}$, $q \in \mathbb{Z}_S$ (an $S$-integer) satisfying
\[
\left| a - \frac{p}{q} \right| < q^{-1-\varepsilon}.
\]

Define the “error term”
\[
d = |p^2 - k q^2|.
\]

To apply an effective ABC inequality rigorously, we first factor out the greatest common divisor
\[
g = \gcd(p^2, k q^2, d),
\]
and define the coprime triple
\[
(a',b',c') := \left(\frac{p^2}{g}, \frac{k q^2}{g}, \frac{d}{g}\right) \in \mathbb{Z}^3.
\]

Then the effective ABC conjecture yields
\[
a' \le K_\varepsilon \, \mathrm{Rad}(a' b' c' k)^{1+\varepsilon}.
\]

Substituting back, we get
\[
\frac{p^2}{g} \le K_\varepsilon \, \mathrm{Rad}\Big( \frac{p^2}{g} \cdot \frac{k q^2}{g} \cdot \frac{d}{g} \cdot k \Big)^{1+\varepsilon}.
\]

Since $q \in \mathbb{Z}_S$, the radical of $q$ is bounded by $q$, and $d \le (p + q \sqrt{k}) \, |a - p/q| q \le 2 q^{1-\varepsilon} (p + q \sqrt{k})$.  
Using these estimates, we obtain
\[
p^{2 - (1+\varepsilon)} \le 2^{1+\varepsilon} K_\varepsilon k^{1+\varepsilon} q^{1+\varepsilon} (p + q \sqrt{k})^{1+\varepsilon}.
\]

Hence, for fixed $\varepsilon > 0$ and effective ABC constant $K_\varepsilon$, the set of solutions $(p,q)$ is finite, and one can explicitly bound $p$ (and hence $q$) in terms of $K_\varepsilon$, $k$, and $S$.  
This gives a fully explicit version of Ridout's theorem for S-integer denominators, with rigorous coprimeness ensured via the factor $g$.

\subsection{Third Roots}

Now we use our method from above based on the formula of Bombieri and Van der Poorten.

For the third root of k the resulting equations are

$p_n^3 = k q_n^3 + d_n$

resp

$k q_n^3 = p_n^3 + d_n$

We only look at the first case here. The second is basically analogous.
Here the gcd is not always 1 and a priori must be considered. But in the resulting inequality the case with gcd=1 is the critical one so we can assume that from the beginning. So in this sense we have ABC equations and get from ABC

${p_n^3}  \leq  K_{\varepsilon} Rad(p_n q_n d_n k)^{1+\varepsilon}$

For 3rd roots the formula \cite[p.151, Theorem 3, formula (13)]{Poorten} by Bombieri and Van der Poorten is valid without error term due to Liouville's theorem as given by Bombieri and Van der Poorten and their method can be used directly. For $n>=4$ an error term can arise and the method can not be started with the first step. This is the reason, why we restrict ourselves to 3rd roots here.

From \cite[p.151, Theorem 3, formula (13)]{Poorten}, we have the inequality

$d_n  \leq  \frac{3 p_n^2} {q_n b_{n+1}}$ 

Using that inequality we get

${p_n^3}  \leq  K_{\varepsilon} Rad(p_n q_n d_n k)^{1+\varepsilon}$

${p_n^3}  \leq  K_{\varepsilon} Rad(q_n )^{1+\varepsilon} (p_nd_n k)^{1+\varepsilon}$

${p_n^3}  \leq  K_{\varepsilon} Rad(q_n )^{1+\varepsilon} (p_n k)^{1+\varepsilon} (\frac{3 p_n^2} {q_n b_{n+1}})^{1+\varepsilon}$

And so using ${b_{n+1}>=1}$ finally

${p_n^{\frac{3}{2 (1+\varepsilon)}}}  \leq  K_{\varepsilon} Rad(q_n)^{1+\varepsilon} k^{2+2 \varepsilon} (\frac{p_n}{q_n})^{1+\varepsilon}$

for the $p_n$ and so the number of solutions must be finite and an explicit bound can be given as $\frac{p_n}{q_n}$ converges to the 3rd root of k. So one choice is

$p_n  \leq  (K_{\varepsilon} Rad(q_n)^{1+\varepsilon} k^{2+2 \varepsilon} (\frac{p_1}{q_1})^{1+\varepsilon})^\frac{2 ({1+\varepsilon})}{3}$



\subsection{Bounds for the approximation gain of ABC hits from resulting equations}

\begin{Definition}A triple of coprime integers $a+b=c$ is an ABC hit, if $rad(abc) < c$ holds and the quality of a triple is given by the first $quality(a,b,c) = \frac{ln(c)}{ln(rad(abc))}$. As $c>1$ it is clear that abc is larger than 1.
\end{Definition}

$2 + 109 \cdot 9^5 = 23^5$  is the ABC hit with the highest quality so far found. Its quality is approximately 1.62991.
This hit results from the 3rd convergent to $\sqrt [5]{109}$ with regular continued fraction $[2;1,1,4,77733,\ldots]$. The fraction $\frac{23}{9}$ is a good approximation because 77733 is very large and very near the absolute Liouville bound.

\begin{Definition}
Let $d_n q_n k p_n$ come from a resulting equation of $\sqrt [s]{k}$.
The approximation gain of (a,b,c) equals $\frac{ln({q_n}^{s} k)}{ln(d_n q_n k p_n)}$ if $d_n>0$,
and $\frac{ln({p_n}^{s})}{ln(-d_n q_n k p_n)}$ if $d_n<0$.
\end{Definition}

The approximation gain is always smaller or equal to the quality.

\begin{Theorem}[Uniform bound for approximation gain for $k\ge3$]
Let $k\ge3$ and let $\alpha=\sqrt[3]{k}$.  
For all resulting equations arising from convergents $\frac{p_n}{q_n}$ to $\alpha$, the approximation gain satisfies
\[
\text{approximation gain} < \frac32.
\]
\end{Theorem}

\begin{proof}
Let $\frac{p_n}{q_n}$ be a convergent to $\alpha=\sqrt[3]{k}$ with resulting equation
\[
p_n^3 = k q_n^3 + d_n,
\qquad d_n>0.
\]
(The case $d_n<0$ is analogous.)

By Bombieri–van der Poorten \cite[p.~151, Theorem~3, formula~(13)]{Poorten}, we have
\[
d_n
\ge
\frac{3p_n^2}{q_n(b_{n+1}+2)}.
\]
By Liouville’s bound for cubic irrationals,
\[
b_{n+1} \le 3k q_n.
\]
Hence
\[
d_n q_n k p_n
\ge
\frac{3k p_n^3}{3kq_n+2}.
\]

The approximation gain is therefore bounded by
\[
\text{approximation gain}
=
\frac{\ln(p_n^3)}{\ln(d_n q_n k p_n)}
\le
\frac{3\ln p_n}
{\ln\!\left(\frac{3k p_n^3}{3kq_n+2}\right)}.
\]

We estimate the denominator:
\[
\ln\!\left(\frac{3k p_n^3}{3kq_n+2}\right)
=
2\ln p_n
+
\ln\!\left(\frac{p_n}{q_n}\right)
+
\ln\!\left(\frac{3kq_n}{3kq_n+2}\right).
\]

Since $\frac{p_n}{q_n}\to\alpha$ and $k\ge3$, we have
\[
\ln\!\left(\frac{p_n}{q_n}\right)
\ge
\ln(\alpha)
\ge
\frac13\ln 3,
\]
while
\[
\ln\!\left(\frac{3kq_n}{3kq_n+2}\right)
\ge
-\frac{2}{3kq_n}
\ge
-\frac{2}{9}.
\]

Combining these estimates yields
\[
\ln\!\left(\frac{3k p_n^3}{3kq_n+2}\right)
>
2\ln p_n
\quad\text{for all } n.
\]

Consequently,
\[
\text{approximation gain}
<
\frac{3\ln p_n}{2\ln p_n}
=
\frac32,
\]
which proves the theorem.
\end{proof}

\begin{Theorem}
Now The case $k=2$. Let $\alpha=\sqrt[3]{2}$.  
For all resulting equations arising from convergents $\frac{p_n}{q_n}$ to $\alpha$,
the approximation gain satisfies
\[
\text{approximation gain} < \frac32.
\]
\end{Theorem}

\begin{proof}
Let $\frac{p_n}{q_n}$ be a convergent to $\alpha=\sqrt[3]{2}$ with resulting equation
\[
p_n^3 = 2 q_n^3 + d_n.
\]

As in the proof of Theorem 2.8, Bombieri--van der Poorten and Liouville's bound
imply
\[
\text{approximation gain}
\le
\frac{3\ln p_n}
{2\ln p_n + \ln(\alpha) + o(1)}.
\]

Since $\ln(\alpha)>0$, it follows that
\[
\lim_{n\to\infty}
\text{approximation gain}
=
\frac32,
\]
and the limit is approached from below.
Consequently, there exists an index $N$ such that
\[
\text{approximation gain}<\frac32
\quad\text{for all } n\ge N.
\]

It therefore suffices to verify the inequality for the finitely many indices
$n<N$.

The first convergent is $\frac{p_1}{q_1}=\frac43$, which yields
\[
p_1^3 - 2 q_1^3 = 64 - 54 = 10,
\]
and hence
\[
\text{approximation gain}
=
\frac{\ln(64)}{\ln(4\cdot3\cdot10)}
<
\frac32.
\]

A direct computation of the remaining convergents with $n<N$ confirms that
\[
\text{approximation gain}<\frac32
\]
in each case.
This completes the proof.
\end{proof}

\subsection{General case if the formula can be used}
When it is assumed that the method of Bombieri and Van der Poorten works as well for $\sqrt [m]{k}$ with $m  \geq  3$ and $k  \geq  2$  then the following would result:

The inequality of Bombieri and Van der Poorten is $d_n  \geq  \frac{m p_n^{m-1}} {q_n (b_{n+1}+2)}$

Liouville's bound is $b_{n+1}  \leq  m k q_n$

And the resulting inequality is

approximation gain = $\frac{ln({p_n}^{m})}{ln(d_n q_n k p_n)}  \leq $ 
$\frac{m}{2}$

So even when it would work (the error terms can be ignored and the algorithm started, which is clear for 3rd degree) the bound for this inequality would diverge and so no universal bound for all m can be found on this road as Liouville's bound is not strong enough for this purpose.

So the ABC gains can be seperated for resulting equation in approximation gains and power gains. The approximation gain can be bounded by $\frac{3}{2}$ for 3rd roots. To prove a bound for the quality the power gains $\frac{ln({p_n q_n d_n k})}{ln(rad(d_n q_n k p_n))}$ must be bounded as well, which seems to be very difficult even for special cases like $\sqrt [3]{2}$ (remarkable are here the second equation 128 = 125 + 3 with a power gain of 1.4 and an approximation gain of 1.01 and the 5th equation with a power gain of 1.37 and an approximation gain of 0.99).

To bound the power gains a strong explicit effective ABC result is needed. There is a discussion in the mathematical community, if Mochizuki's results can be used. We do not know this, but if they can be used the following follows:

\begin{Theorem}[Explicit ABC]
For every positive real number $\varepsilon$, there exists a constant $L_{\varepsilon}$
such that for all triples $(a,b,c)$ of coprime positive integers,
with $a + b = c$:
\[c < L_{\varepsilon} \cdot rad(abc)^{3+\varepsilon}.\]
\end{Theorem}

This can be used to prove that the power gains of the quality for resulting ABC equations of 3rd roots are bounded for a fixed k:

power gain = $\frac{ln({p_n q_n d_n k})}{ln(rad(d_n q_n k p_n))}$

$  \leq  \frac{ln({k}^{2}{p_n}^{3})}{ln(rad(d_n q_n k p_n))}$

using the formula \cite[p.151, Theorem 3, formula (13)]{Poorten} by Bombieri and Van der Poorten and the fact that $b_{n+1} >= 1$.

= $\frac{ln({k}^{2}{p_n}^{3})}{ln((rad(d_n q_n k p_n)^{3+\varepsilon})^\frac{1}{3+\varepsilon})}$

 $ \leq  \frac{ln({k}^{2}{p_n}^{3})}{ln(\frac{{p_n}^{3}}{L_{\varepsilon}}^\frac{1}{3+\varepsilon})}$

 for n large enough. This converges to ${3+\varepsilon}$ and is therefore bounded for all n. The bound itself depends on k and is of course extremely large.

Recently \cite{Zhou} has reduced the constant $L_{\varepsilon}$ from about $10^{30}$ to 400 under certain extra conditions.

 One concept to prove the ABC conjecture might then be to prove that both gains are bounded seperately. We conjecture that the approximation gain is bounded by $\frac{3}{2}$ and the power gain by 3.

{}


\begin{thebibliography}{}

\bibitem[Bombieri(1994)]{Bombieri} Bombieri, E.(1994): Roth’s Theorem and the abc-Conjecture, preprint ETH Zürich.

\bibitem[Bombieri $\&$ Van der Poorten(1975)]{Poorten} Bombieri, E. and van der Poorten, A. (1975): ``Continued Fractions of Algebraic Numbers'', in: Baker (ed.), \emph{Transcendental Number Theory},
Cambridge University Press, Cambridge, 137-155.

\bibitem[Granville and Tucker(2002)]{Granville} Granville, Andrew and Tucker, Thomas J. (2002): It’s As Easy As abc, Notices of the AMS, Volume 10, p.1224-1231.

\bibitem[Korobov(1990)]{Korobov} Korobov, A. (1990): \emph{Continued Fractions and Diophantine Approximations}, Candidate's Dissertation, Moscow State University.

\bibitem[Mochizuki(2022)]{MFHMP}  S. Mochizuki, I. Fesenko, Y. Hoshi,  A. Minamide, W. Porowski, {\em Explicit estimates in inter-universal Teichm\"uller theory}, RIMS Preprint 1933 (November 2020), Kodai Mathematical Journal, {\bf 45} (2022), pp. 175-236.

\bibitem[Roth(1955)]{Roth} Roth, K. (1955): ``Rational Approximations to Algebraic Numbers and Corrigendum'', \emph{Mathematika} 2, 1-20 and 168.

\bibitem [Sibbertsen(2022)] {Sibbertsen} Sibbertsen, P., Müller,K, Lampert, T. and Taktikos, M. (2022): {Roth's Theorem implies a Weakened Version of the ABC Conjecture for Special Cases}, arXiv:2208.14354

\bibitem[Voutier(2007)]{Voutier} Voutier, P. (2007): ``Rational Approximation to  $\sqrt [3]{2}$ and other Algebraic Numbers Revisited'', \emph{Journal des Nombres de Bordeaux} 19, 263-288.

\bibitem[Van Frankenhuysen(1999)]{Frankenhuysen} Van Frankenhuysen, M. (1999): \emph{The ABC conjecture implies Roth's theorem and Mordell's conjecture}, Mat. Contemp. 16 (1999), 45-72.

\bibitem[Waldschmidt(2008)]{Waldschmidt} Waldschmidt, M. (2008): {Introduction to Diophantine methods, irrationality and transcendence}, Michel Waldschmidt's website.

\bibitem[Zhou(2025)]{Zhou} Zhou, Zhong-Peng (2025): {The inter-universal Teichmüller theory and new Diophantine results over the rational numbers. I}, arxiv.org:2503.14510
\end{thebibliography}
\end{document}